\documentclass{amsart}

\usepackage{amsmath,amssymb,amsfonts,enumerate,amsthm, amscd}

\newcommand{\Ext}{\mbox{Ext}\,}

\newcommand{\pd}{\mbox{pd}\,}
\newcommand{\id}{\mbox{id}\,}
\newcommand{\fd}{\mbox{fd}\,}

\newcommand{\wdim}{\mbox{wdim}\,}
\newcommand{\gldim}{\mbox{gldim}\,}
\newcommand{\gwd}{\mbox{Gwdim}\,}
\newcommand{\ggldim}{\mbox{Ggldim}\,}
\newcommand{\wggldim}{\mbox{wGgldim}\,}
\newcommand{\gpd}{\mbox{Gpd}\,}
\newcommand{\gid}{\mbox{Gid}\,}
\newcommand{\gfd}{\mbox{Gfd}\,}

\newtheorem{theorem}{Theorem}[section]
\newtheorem{lemma}[theorem]{Lemma}
\newtheorem{proposition}[theorem]{Proposition}

\theoremstyle{definition}
\newtheorem{definition}[theorem]{Definition}
\newtheorem{definitions}[theorem]{Definitions}
\theoremstyle{remark}
\newtheorem{remark}[theorem]{Remark}
\newtheorem{remarks}[theorem]{Remarks}

\theoremstyle{Definition and Notation}

\begin{document}
\bibliographystyle{amsplain}

\title[Gorenstein Von Neumann regular rings]{Gorenstein Von Neumann regular rings}

\author{Najib Mahdou}
\address{Najib Mahdou\\Department of Mathematics, Faculty of Science and Technology of Fez, Box 2202, University S.M. Ben Abdellah Fez, Morocco.}

\author{Mohammed Tamekkante}
\address{Mohammed Tamekkante\\Department of Mathematics, Faculty of Science and Technology of Fez, Box 2202, University S.M. Ben Abdellah Fez, Morocco.}

\author{Siamak Yassemi}
\address{Siamak Yassemi\\Department of Mathematics,
University of Tehran, Tehran, Iran\\ and School of Mathematics,
Institute for research in fundamental sciences (IPM), P.~O.~Box
19395-5746, Tehran, Iran.}




\begin{abstract}

In this paper, we study the rings with zero Gorenstein weak
 dimensions, which we call them Gorenstein Von Neumann regular rings.

\end{abstract}

\maketitle

\section{Introduction}
Throughout this paper, all rings are commutative with identity
element, and all modules are unital.

Let $R$ be a ring, and let $M$ be an $R$-module. As usual we use
$\pd_R(M)$, $\id_R(M)$ and $\fd_R(M)$ to denote, respectively, the
classical projective, injective and flat dimensions of $M$. By
$\gldim(R)$ and $\wdim(R)$ we denote, respectively, the classical
global and weak dimensions of R. It is by now a well-established
fact that even if R is non-Noetherian, there exists Gorenstein
projective, injective and flat dimensions of M, which are usually
denoted by $\gpd_R(M)$,
 $\gid_R(M)$ and $\gfd_R(M)$, respectively. Some references are
 \cite{Bennis and Mahdou1, Bennis and Mahdou2, Christensen, Christensen
 Frankild, Enochs, Enochs2, Eno Jenda Torrecillas, Holm}.\\

 Recently in \cite{Bennis and Mahdou2}, the authors started the study of
global Gorenstein projective (resp. global Gorenstein injective,
global weak) dimensions of $R$, denoted by $\gpd(R)$ (resp.
$\gid(R)$, $\gwd(R)$) and defined as the supremum of Gorenstein
projective (resp. Gorenstein injective, Gorenstein flat) dimensions
of $M$ where $M$ ranges over all $R$-modules.

It is known that for any ring R, $\gwd(R)\leq \gid(R) = \gpd(R)$,
see \cite[Theorems 1.1 and Corollary 1.2(1)]{Bennis and Mahdou2}.
So, according to the terminology of the classical theory of
homological dimensions of rings, the common value of $\gpd(R)$ and
$\gid(R)$ is called Gorenstein global dimension of $R$, and denoted
by $G.gldim(R)$. The Gorenstein global (resp. weak) dimension is
refinement of the classical global (resp. weak) dimension of rings.
That is $\ggldim(R) \leq \gldim(R)$ (resp. $\gwd(R)\leq \wdim(R)$)
with equality holds if $\gldim(R)$ (resp. $\wdim(R)$) is finite, see
\cite[Corollary 1.2(2 and 3)]{Bennis and Mahdou2}).

The rings with zero global Gorenstein dimension (i.e. every module
is Gorenstein projective module) are studied in \cite{Ouarghi},
where the authors called them Gorenstein semisimple.

In this paper, motivating by the work in \cite{Ouarghi}, we study
the rings with zero Gorenstein weak dimension. Analogy to the
classical ones we call them Gorenstein Von Neumann regular rings.

It is known result that in a Von Neumann regular ring, every finitely generated projective submodule of a projective $R$-module $P$ is a direct summand of $P$.
One of the result in this paper is an analog of this classical one.
It is shown that $R$ is Gorenstein Von Neumann regular ring if and only if  $R$ is Gorenstein semihereditary ring
that every finitely generated projective submodule of a projective $R$-module $P$ is a direct summand of $P$ (recall that a ring $R$ is called to be Gorenstein semihereditary if it is coherent and
every submodule of a projective module is Gorenstein projective).

It is known that a semisimple ring is exactly a Noetherian Von Neumann
regular ring. The second main result of this paper is an analog of
this classical one. We show that $R$ is Gorenstein semisimple ring if and only if $R$ is Noetherian Gorenstein Von Neumann regular.

Above we have only mentioned the Gorenstein von Neumann regular
rings. One can also define the strongly Gorenstein Von Neumann
regular rings by using the notions of strongly Gorenstein global
dimension, cf. \cite{Bennis and Mahdou1}. All the results concerning Gorenstein Von Neumann regular
rings, have a strongly Gorenstein Von Neumann regular part. We do
not state or prove these "strongly" Gorenstein Von Neumann regular
rings results. This is left to the reader.


\section{Main Results}


\begin{definitions}

A ring $R$ is called  Gorenstein Von Neumann regular ($G$-VNR for
short) if every $R$-module is $G$-flat (i.e., $G-wdim(R) = 0$).

\end{definitions}


\begin{remarks}\label{remarks G-semisimple--G-RVN}\
\begin{enumerate}

\item From \cite[Corollary 1.2(1)]{Bennis and Mahdou2}, we can see that every
$G$-semisimple ring is $G$-Von Neumann regular with equivalence if
the ring is Noetherian (\cite[Theorem 1.2.3.1]{Enochs3}).
\item Every Von Neumann regular ring is $G$-Von Neumann regular
with equivalence if $wdim(-)$ is finite (\cite[Corollary
1.2]{Bennis and Mahdou2}).
\end{enumerate}
\end{remarks}

\begin{definition}(\cite{Stenstrom} and \cite{Garkusha}) Let $R$ be a
ring and $M$ an  $R$-module.
\begin{enumerate}
    \item We say that $M$ is $FP$-injective (or absolutely pure) if $\Ext^1_R(P,M) = 0$ for every finitely
     presented $R$-module $P$.
    \item $R$ is said to be a weakly quasi-Frobenius ring (or also $FC$-ring) if it is coherent and it is self-$FP$-injective
    (i.e., $R$ is $FP$-injective as an $R$-module).
\end{enumerate}
\end{definition}

$R$ is called IF-ring if every injective R-module is flat
\cite{Colby}. Over a commutative rings, the IF-rings and the weakly
quasi-Frobenius rings are the same (\cite[Proposition
4.2]{Stenstrom}).

\begin{lemma}
The following are equivalent:
\begin{enumerate}
    \item Every injective $R$-module is flat (ie., IF-ring).
    \item Every Gorenstein injective $R$-module is Gorenstein
    flat (ie., GIF ring).
\end{enumerate}
\end{lemma}

\begin{proof}
$1\Rightarrow 2.$ Consider an arbitrary  complete injective
resolution:
$$\mathbf{I}: ...\rightarrow I_1\rightarrow I_0\rightarrow I^0\rightarrow
I^1\rightarrow ...$$ Since every injective module is flat,
$\mathbf{I}$ is also an exact flat resolution. Moreover, for every
injective $J$, $\mathbf{I}\otimes_RJ$ still exact (since $J$ is also
flat). Thus, $\mathbf{I}$ is a complete flat resolution. Then, every
Gorenstein injective
$R$-module is Gorenstein flat, as desired.\\
$2\Rightarrow 1.$ Let $I$ be an arbitrary injective $R$-module (and
so Gorenstein injective). Then, by hypothesis, $I$ is Gorenstein
flat. Thus, it can be embedded in a flat module $F$. So, $I$ is a
direct summand of $F$ and so flat, as desired.
\end{proof}

\begin{proposition}
Let $R$ be a coherent commutative ring. The following are
equivalent:
\begin{enumerate}
    \item $\wggldim(R)\leq n$.
    \item $\fd_R(I)\leq n$ for every injective $R$-module $I$ (ie, n-IF).
    \item $\gfd_R(I)\leq n$ for every Gorenstein injective
    $R$-module $I$ (ie., $n$-GIF).
\end{enumerate}
\end{proposition}

\begin{proof} The equivalence $(1\Leftrightarrow 2)$
follows by combining the equivalence \cite[Theorem
7$(1\Leftrightarrow 2)$]{Ding} and the equality \cite[Theorem
3.7(1=2)]{Chen}.\\
$1\Rightarrow 3.$ Obvious.\\
 $3\Rightarrow 2.$ Let $I$ be an arbitrary injective $R$-module.
 By hypothesi, $\gfd_R(I)\leq n$. Thus, from \cite[Theorem 2.10]{Holm}, there is an exact sequence
$0\rightarrow K \rightarrow G \rightarrow I\rightarrow 0$ where
 $G$ is Gorenstein flat and $\fd_R(K)\leq n-1$. For the $R$ module
 $G$ we can pick an short exact sequence $0\rightarrow G
 \rightarrow F \rightarrow G'\rightarrow 0$ where $F$ is flat and
 $G'$ is Gorenstein flat. We have  the following
pushout diagram:
$$\begin{array}{ccccccc}
   & 0 &  & 0 &  &  &  \\
   & \downarrow &  & \downarrow &  & &  \\
   & K& = & K &  &  &  \\
   & \downarrow &  & \downarrow &  &  &  \\
  0\rightarrow & G & \rightarrow & F  & \rightarrow & G' & \rightarrow 0 \\
   & \downarrow &  & \downarrow &  & \parallel &  \\
  0\rightarrow & I & \rightarrow & D & \rightarrow  & G'& \rightarrow 0 \\
 & \downarrow &  & \downarrow &  & &  \\
 & 0 &  & 0 &  &  &  \\
\end{array}$$
Clearly, $\fd_R(D)\leq n$. On the other hand, $I$ is a direct
summand of $D$. Thus, $\fd_R(I)\leq n$ as desired.
\end{proof}


\begin{theorem}\label{G-RVN} The following conditions are equivalent:
\begin{enumerate}
    \item $R$ is $G$-Von Neumann regular.
    \item Every finitely presented $R$-module is $G$-flat.
    \item Every finitely presented $R$-module is $G$-projective.
    \item  Every finitely presented $R$-module embeds in a flat
    $R$-module.
    \item  $R$ is coherent and self-FP-injective ring (i.e., weakly
quasi-Frobenius).
    \item Every injective $R$-module is flat (i.e., $R$ is IF-ring).
\end{enumerate}
\end{theorem}

\begin{proof}  The equivalences $(1\Leftrightarrow 2\Leftrightarrow 3)$ is an
other way to see   \cite[Theorem
6]{Chen}. \\
$1\Rightarrow 4$. Easy by  definition of a Gorenstein
flat modules.\\
$4\Leftrightarrow 5 \Leftrightarrow 6$.  By
\cite[Proposition 2.5 and Theorem 2.8]{Garkusha}.\\
$4\Rightarrow 1$. For any $R$-module $M$, assemble any flat
resolution of $M$ with its any injective resolution, we get an
exact sequence of flat $R$-modules (by hypothesis), which is also
by hypothesis a complete flat resolution. This means that $M$ is
Gorenstein flat.
\end{proof}


\begin{remark} From Theorem \ref{G-RVN}, we deduce that a $G$-Von
Neumann regular ring is always coherent.
\end{remark}


\begin{proposition}\label{caracte 2 GRVN Coherent}
If $R$ is a coherent ring. Then, the following conditions are
equivalent:
\begin{enumerate}
    \item $R$ is a $G$-Von Neumann regular ring.
    \item For every finitely generated ideal $I$, the module $R/I$
    is $G$-flat.
 \item Every injective $R$-module is flat.
 \item Every flat $R$-module is FP-injective.
  \item Every $FP$-injective $R$-module is flat.
\end{enumerate}
\end{proposition}

\begin{proof}
$1\Rightarrow 2$. Easy.\\
$2\Rightarrow 3.$ Follows from \cite[Theorem 3.4]{Holm} and
\cite[Theorem 1.3.8]{Glaz}. \\
$3\Rightarrow 4 \Rightarrow 5 \Rightarrow 1$. Follows from
\cite[Theorem 3.5 and 3.8]{Ding} and Theorem \ref{G-RVN}.
\end{proof}



Recall that a ring $R$ is called to be Gorenstein semihereditary if
it is coherent and every submodule of a projective module is
Gorenstein projective. In other words, the ring $R$ is Gorenstein
semihereditary if $R$ is coherent and $G.wdim(R)\leq 1$
(\cite{MahdouTamekkante2}).

The following two properties ``every finitely generated proper ideal
of $R$ has nonzero annihilator" and ``every finitely generated
projective submodule of a projective $R$-module $P$ is a direct
summand of $P$" are equivalent by \cite[Theorem 5.4]{Bass}. For
short, we call such ring a CH-ring. It is known that a Von Neumann
regular ring is exactly a semihereditary ring which is a CH-ring.
The first main result of this paper is an analog of this classical
one.


\begin{theorem}\label{G-semihere+ CH}
The following statements are equivalents:
\begin{enumerate}
    \item $R$ is Gorenstein Von Neumann regular;
    \item $R$ is Gorenstein semihereditary which is CH-ring.
\end{enumerate}
\end{theorem}


The proof of the theorem involves the following Lemma which is a new
characterization of a CH-ring.


\begin{lemma}\label{CH-GCH}
 Let $R$ be a ring. The following conditions are equivalents:
 \begin{enumerate}
    \item Every finitely generated projective submodule of a
    $G$-projective $R$-module is a direct summand.
    \item $R$ is a CH-ring.
 \end{enumerate}
\end{lemma}

\begin{proof}
$1\Rightarrow 2$ Obvious  since every projective module is $SG$-projective and every $SG$-projective module is $G$-projective.\\
$2 \Rightarrow 1$ Assume that $R$ satisfied the $(CH)$-property and
let $P$ be a finitely presented projective
 submodule of a $G$-projective module $G$. We claim that $P$ is a direct summand of $G$. Pick a short exact sequence $0 \rightarrow G
\rightarrow Q \rightarrow G' \rightarrow 0$ where $Q$ is a
projective $R$-module and $G'$ is a $G$-projective $R$-module.
Thus, we have the following pull-back diagram with exact rows and
columns:
$$\begin{array}{ccccccc}
  & 0 &  & 0 &  &  &  \\
   & \downarrow &  & \downarrow &  &  &  \\
   & P & = & P &  &  &  \\ & \downarrow &  & \downarrow &  &  &  \\
  0\rightarrow & G & \rightarrow & Q & \rightarrow & G' & \rightarrow 0 \\
   & \downarrow &  & \downarrow &  & \parallel &  \\
 0\rightarrow & G/P & \rightarrow & K & \rightarrow & G' & \rightarrow 0 \\
  & \downarrow &  & \downarrow &  &  &  \\
 &  0&  & 0 &  &  &  \\
\end{array}$$
From the middle vertical sequence we conclude that $P$ is isomorphic
to a direct summand of $Q$ (since $R$ is a CH-ring). So, $K$ is also
isomorphic to a direct summand of $Q$. Consequently, it is
projective. This implies, from \cite[Theorem 2.5]{Holm}, that $G/P$
is $G$-projective and so, $\Ext(G/P,P)=0$. Therefore, the short
exact sequence $0\rightarrow P \rightarrow G  \rightarrow  G/P
\rightarrow 0$ splits, that means that $G=P\oplus G/P$, as desired.
\end{proof}

\noindent\textit{Proof of Theorem \ref{G-semihere+ CH}.} \\
Firstly assume that $R$ is a Gorenstein Von Neumann regular ring.
Clearly, $R$ is Gorenstein semihereditary. So, we have to prove
that $R$ satisfied the $(CH)$-property. Let $P$ be a finitely
generated projective submodule of a projective module $Q$. We
claim that $P$ is a direct summand of $Q$. Let $Q'$ be a
projective $R$-module such that $Q\oplus Q'$ is a free $R$-module.
We have the short exact sequence

$$\begin{array}{ccccccc}
  0 \longrightarrow  & P &\longrightarrow & Q\oplus Q'& \longrightarrow &Q/P \oplus Q'& \longrightarrow 0\\
   & x & \mapsto & (x,0) &  &  &  \\
   &  &  & (x,y) &\mapsto & (\bar{x},y) &  \\
\end{array}$$

 We can identify $P$ to be a
submodule of $L=Q\oplus Q'$. Now, let $L_0$ be a finitely
generated free direct summand of $L$ such that $P\subset L_0$
(that exists since $P$ is finitely generated submodule of $L$) and
let $L_1$ be a free module such that $L=L_0\oplus L_1$ . The
$R$-module $L_0/P$ is finitely presented and so it is
$G$-projective by Theorem \ref{G-RVN}. Therefore, $L_0/P$ is
projective by \cite[Proposition 2.27]{Holm} (since
$pd_R(L_0/P)\leq 1$). Now, we consider  the following pull-back
diagram with exact rows and columns:
$$\begin{array}{ccccccc}
   &  &  & 0 &  & 0 &  \\
   &  &  & \downarrow &  & \downarrow &  \\
  0\rightarrow & P & \rightarrow & L_0 & \rightarrow & L_0/P & \rightarrow 0 \\
   & \| &  & \downarrow &  & \downarrow &  \\
 0\rightarrow & P & \rightarrow & L & \rightarrow & Q/P\oplus Q' & \rightarrow 0 \\
   &  &  & \downarrow &  & \downarrow &  \\
   &  &  & L_1 & = & L_1 &  \\
   &  &  & \downarrow &  & \downarrow &  \\
      &  &  & 0 &  & 0 &  \\
\end{array}$$
 From the right  exact sequence we deduce that $Q/P\oplus Q'$ is
 projective (since $L_0/P$ and $L_1$ are projective). Then, $Q/P$
 is projective and so $Q=Q/P\oplus P$, as desired.\\
Conversely, assume  that $R$ is a $G$-semihereditary ring which
satisfied the $(CH)$-property. From Theorem  \ref{G-RVN}, to prove
that $R$ is $G$-Von Neumann regular, we have to prove that every
finitely presented $R$-module is $G$-projective. So, let $M$ be a
finitely presented $R$-module and pick a short exact sequence of
$R$-modules $0 \rightarrow K \rightarrow P \rightarrow
M\rightarrow 0$, where $P$ is a finitely generated projective
$R$-module and $K$ is a finitely generated $R$-module. By
\cite[Theorem 7]{Chen}, the $R$-module $K$ is $G$-projective
(since $R$ is $G$-semihereditary then is also coherent). Then,
there is an exact sequence of $R$-modules $0 \rightarrow K
\rightarrow Q \rightarrow K' \rightarrow 0$ where $Q$ is
projective and $K'$ is $G$-projective. Now, consider the pull-bach
diagram with exact rows and columns:
\begin{center}
$\begin{array}{ccccccc}
 & 0 &  &0&  &    &  \\
    & \downarrow &  & \downarrow &  &   &  \\
  0\rightarrow &K & \rightarrow & P & \rightarrow & M & \rightarrow 0 \\
  & \downarrow &  &\downarrow &  & \|   &  \\
 0\rightarrow &Q & \rightarrow &Z & \rightarrow & M & \rightarrow 0 \\
 & \downarrow &  &\downarrow &  &    &  \\
  & K' & = & K'&  &  &  \\
   & \downarrow &  &\downarrow &  &    &  \\
    & 0 &  &0&  &    &  \\
\end{array}$
\end{center}
The middle vertical short exact sequence $0 \rightarrow P
\rightarrow Z \rightarrow K'\rightarrow 0$ splits since $K'$ is
$G$-projective and $P$ is projective (by \cite[Theorem
2.20]{Holm}). Then $Z$ is $G$-projective  ($Z= P\oplus K$). Now,
from  Lemma \ref{CH-GCH}, $Q$ is isomorphic to a direct summand of
$Z$ (since $R$ satisfied the $(CH)$-property).
 Then, $M$ is also isomorphic to a direct summand of $Z$ and so $M$ is $G$-projective by \cite[Theorem 2.5]{Holm}, as
 desired.\\ The strongly special cases is from the fact that the
 $G$-semihereditary ring (resp. the $G$-Von Neumann regular ring)
 is $SG$-semihereditary ring (resp. the $SG$-Von Neumann regular
 ring) if, and only if, every $G$-flat module is $SG$-flat. $\Box$

We know that a semisimple ring is exactly a Noetherian Von Neumann
regular ring. The second main result of this paper is an analog of
this classical one.

 \begin{proposition}
 Let $R$ be a ring. The following statement are equivalents.
 \begin{enumerate}
    \item $R$ is Gorenstein semisimple;
    \item $R$ is Noetherian Gorenstein Von Neumann
    regular.
 \end{enumerate}
 \end{proposition}

 \begin{proof} From \cite[Proposition 2.6]{Bennis and
 Mahdou2}, the quasi Frobenius rings and the Gorenstein semisimple rings are the
 same; then, they are Noetherian.  On the other hand, by
 \cite[Theorem 12.3.1]{Enochs3}, if $R$ is Noetherian we have
 $G.wdim(R)=G.gldim(R)$. Thus, $R$ is Gorenstein
 semisimple if, and only if, $R$ is a Noetherian  Gorenstein Von Neumann
    regular.\\
    So, to finish the proof we have to prove the strongly
    particular cases.\\
 Assume that $R$ is Gorenstein
 semisimple and  we claim that $R$ is Gorenstein Von Neumann
 regular. Let $M$ be an arbitrary $R$-module. We claim that $M$ is Gorenstein flat. From above, $R$ is
 Gorenstein Von Neumann regular and so $M$ is Gorenstein flat; hence, for any injective module $I$, $Tor(M,I)=0$. On the other hand
 since $R$ is Gorenstein semisimple, $M$ is
 Gorenstein projective and so, from \cite[Proposition 2.9]{Bennis and Mahdou1}, there is an exact sequence
 $0\rightarrow M \rightarrow P \rightarrow M  \rightarrow 0$ where
 $P$ is projective (then flat). Consequently, by \cite[Proposition 3.6]{Bennis and
 Mahdou1}, $M$ is Gorenstein flat, as desired.\\
 Conversely, assume that $M$ is a Noetherian Gorenstein Von Neumann regular ring and let  $M$
 an arbitrary $R$-module. We claim that $M$ is Gorenstein
 projective. From the first part of the proof, $R$ is Gorenstein
 semisimple ring. Thus, $M$ is Gorenstein projective and so for
 any projective module $P$ we have $Ext(M,P)=0$. On the other
 hand, $M$ is Gorenstein flat module and so, from \cite[Proposition 3.6]{Bennis and Mahdou1}, there is an
 exact sequence  $0\rightarrow M \rightarrow F \rightarrow M
 \rightarrow 0$ where $F$ is flat. Using \cite[Corollary 2.7]{Bennis and
 Mahdou2} and the fact that $R$ is Gorenstein semisimple, $F$ is
 also projective. Hence, $M$ is Gorenstein projective
 module (\cite[Proposition 2.9]{Bennis and Mahdou1}), as desired.
 \end{proof}

The final result of this paper is an analog of the classical one
which say that a Von Neumann regular domain is a field. \\


 \begin{proposition}Every Gorenstein Von Neumann regular  domain is a field.
\end{proposition}

\begin{proof} Assume that $R$ is a Gorenstein Von Neumann regular domain. We
claim that $R$ is a field. So, let $x$ be a nonzero element of
$R$. We have to prove that $x$ is invertible.
     Clearly, $\pd(R/xR)\leq 1$ since $R$ is a domain. On the
     other hand, from Proposition \ref{G-RVN}, $R/xR$ is $G$-projective since it is finitely presented and
     $R$ is Gorenstein Von Neumann regular. Therefore,
     from \cite[Proposition 2.27]{Holm}, $R/xR$ is projective. Thus,  the  short exact sequence
     $0\rightarrow xR\rightarrow R\rightarrow R/xR\rightarrow 0$ splits and so $R\cong xR\oplus
     R/xR$. Localizing this isomorphism with  a maximal ideal $\mathcal{M}$ of $R$, we
     obtain
$R_{\mathcal{M}}\cong
     (x/1)R_{\mathcal{M}}\oplus (R/xR)_{\mathcal{M}}$.  But,  $(x/1)R_{\mathcal{M}}\neq \{0\}$  since $R$ is a domain. So,  $(R/xR)_{\mathcal{M}}=\{0\}$
    for every maximal ideal $\mathcal{M}$ since $R_{\mathcal{M}}$ is local. This implies that   $
   R/xR=0$ which means that $x$ is invertible, as desired.
\end{proof}


\begin{remark} A local $G$-Von Neumann regular  (or  $G$-semisimple) ring is not necessarily a
   field. For example, we can take $R=K[X]/(X^2)$ which is a local
   $G$-semisimple (then $G$-Von Neumann regular) ring but not a field.
   \end{remark}



\end{document}